\documentclass[12pt]{article}

\usepackage{amsmath,amsbsy,amscd,amssymb,graphicx, epsfig,color}
\usepackage[english]{babel}
\usepackage{graphics}
\usepackage[left=1.5cm,top=2.5cm,right=1.5cm,bottom=2.5cm]{geometry}

\def \sin {\hbox{\rm  sin }}
\def \cos {\hbox {\rm  cos }}

\newcommand\N{\mathbb{N}}

\newcommand\R{\mathbb{R}}

\newtheorem{Def}{{\mbox{$\;\;\;\;\;\,$}}Definition}[section]
\newtheorem{Th}{{\mbox{$\;\;\;\;\;\,$}}Theorem}[section]

\newtheorem{Lem}{{\mbox{$\;\;\;\;\;\,$}}Lemma}[section]

\newtheorem{Ex}{{\mbox{$\;\;\;\;\;\,$}}Example}[section]
\newtheorem{Rem}{{\mbox{$\;\;\;\;\;\,$}}Remark}[section]

\newcommand{\qed}{\nobreak \ifvmode \relax \else
     \ifdim\lastskip<1.5em \hskip-\lastskip
      \hskip1.5em plus0em minus0.5em \fi \nobreak
      \vrule height0.75em width0.5em depth0.25em\fi \bigskip}

\newenvironment{proof}{{\noindent\it\underline{Proof}}:}{\hfill{$\blacksquare$}}

\title{On almost periodic solutions for a model of hematopoiesis  with an oscillatory circulation loss rate. }

\author{\small{Roc\'io Balderrama\footnote{Fax:  +54-11-45763335.\vskip0.5pt  
			E-mail addresses: rbalde@dm.uba.ar (R. Balderrama)
			\vskip0.1pt 
			This work was fully supported by the project UBACyT 20020120100029BA}}}
\providecommand{\keywords}[1]{\small{\textbf{{Keywords---}}} #1}

\begin{document}
	\date{}
	\maketitle
	
{\centerline{\small\textit{Departamento de Matem\'atica, Facultad de Ciencias Exactas y Naturales}}}

{\centerline{\small\textit{Universidad de Buenos Aires \& IMAS-CONICET - Argentina.}}}

	\begin{abstract}
		We establish and prove a  fixed point theorem from which some sufficient conditions on the existence of positive almost periodic
		solutions for a model of hematopoiesis with oscillatory circulation loss rate are deduced. 
		Some particular assumption under the nonlinearity of the equation has been previously considered by authors as fundamental for the study of almost periodic solutions of the model. The aim of this paper is to establish results without such assumption. Some examples are given to illustrate our results.
	\end{abstract}
	
		

{\thispagestyle{empty}} 

\section{Introduction}
\label{intro}

The following nonlinear model was proposed
by Mackey and Glass \cite{MG}  
to describe the regulation of hematopoiesis: 
\begin{equation}
\label{MG} 
x'(t) = \sum_{k=1}^M r_k(t)\frac{x^m(t-\tau_k(t))}{1+x^n(t-\tau_k(t))}-b(t)x(t).
\end{equation}
where $x(t)$ is the concentration
of cells in the circulating blood,  $f(x(t-\tau_k(t)))=r_k(t)\frac{x(t-\tau_k(t))^m}{1+x(t-\tau_k(t))^n}$  
is the flux function of cells into the blood stream and the delay $\tau_k(t)$ is the time  between the start of cellular production in the 
bone marrow and the release of mature cells into the blood stream at time t.

The existence of almost periodic solutions for (\ref{MG}) has been extensively studied by a number of authors, for example, see \cite{AB,DLN,DZ,Jiang,L} and references therein. In these works, we find that authors assume the following condition:

$\hspace{0.1mm}(H_0) \qquad  0\le m \le 1,$ 

\noindent
as fundamental for the study of  (\ref{MG}). In addition, in \cite{DLN,DZ, L}, authors proposed the \textit{open problem}  of extending the existence results for (\ref{MG}) to the case $m>1$.

As pointed out in   \cite{BB2, Jiang}, during some seasons the death rate becomes greater than the birth rate, so equations with oscillating coefficients have more realistic significance.
In this work, with respect to equation (\ref{MG})  we will assume that $\lambda_{k}, n, m$ are positive constants, $m>1$,  $r_{k}(t), b(t)$ and $ \tau_{k}(t)\in AP(\R)$, $r_{k}(t)$  are positive, $ \tau_{k}(t)$ is nonnegative and   $b(t)$ is oscillating. 
To the best of our knowledge it is the first time to focus on the dynamic behavior of (\ref{MG}) without condition $(H_0)$.

Throughout the paper,   for a bounded continuous function $f$ we employ de notation 

$$f^+ = \sup_{t\in\R}f(t) \qquad \hbox{ and } \qquad f^- = \inf_{t\in\R}f(t).$$ 
Moreover, it will be assumed 
that 

\begin{equation}
\label{supremopositivo}
M[b] = \lim_{t\to +\infty}\frac{1}{T}\int_t^{t+T}b(s)ds>0,\hspace{0.8cm} \upsilon:=\max_{1\le k\le M}\left\{\tau_k^+\right\}>0 \hspace{0.5cm} \hbox{ and }\hspace{0.5cm}  r_j^->0 \,\,\,\hbox{ for some } j
\end{equation}
and there exist a positive constant $F^s$, such that

\begin{equation}
\label{oscilatoria}
e^{-\int_s^t b(u)du}\le F^s e^{-\int_s^t b^*(u)du}, \,\,\hbox{ for all } t,s \in \R \hbox{ such that } t\geq s,
\end{equation}
where $ b^*:\R \to (0,+\infty)$ is a bounded continuous function with positive infimum.

\begin{Rem}\label{oscilatoria2}
	It is worth noticing that when $b^->0$, inequation (\ref{oscilatoria}) is fulfilled with $ F^s = 1$ and ${b^*}(t) \equiv b(t)$. 
	Thus, our  techniques are applicable when $b(t)$ is not oscillating.  
\end{Rem}

\section{Preliminaries}
\label{prelim}
\begin{Def}(Corduneanu \cite{Co}) 
	\label{bohr}
	Let $X$ be a Banach space. A function $f:\R\to X$  is called \textbf{almost periodic} if for any $\epsilon>0$ there exists a number $l(\epsilon)>0$
	such that any interval on $\R$ of length $l(\epsilon)$ contains at least one point  $\xi$ with
	the property that 	$$||f(t+\xi)-f(t)||<\epsilon \qquad\,\hbox{for all}\; t\in\R.$$ 
\end{Def}

\begin{Def}
	Let $X$ be a real Banach space.  
	A nonempty closed set $C\subset X$ is called a \textbf{cone} if the following conditions are fulfilled:
	$$\hbox{ (a) }  C+C\subset C \qquad\hbox{ (b) } C\cap -C =\{0\} \qquad\hbox{ (c) } C \hbox{ is convex, }
	$$
	where  $0$ denotes the zero element of $X$.
	
	Every cone
	$C$  induces a partial order  $\le$ in $X$ given  by 	
	$$x\le y \hbox{ if and only if } y-x\in C.$$
	If $x\le y $ and $x\neq y$, we write $x<y$. A set $\{z\in X/ x\le z \le y\}$ is called an \textbf{ order interval} and shall be denoted as $[x, y]$. 
	The interior of $C$ shall be denoted by $C^{\circ}$. 
	A cone $C$ is called \textbf{normal} if there exists a constant $N>0$ such that
	$$
	0\le x\le y \hbox{ implies that } ||x||\le N||y||.
	$$
	The smaller constant $N$ satisfying the inequality is called the \textbf{normal constant} of $C$.
\end{Def}

We denote by $AP(\R)$ the Banach space of almost periodic real functions defined on $\R$,
equipped with the usual supremum  norm.  Also, we denote 

$$P:=\{x\in AP(\R) : x(t)\geq 0, \forall t\in \R\},$$
the normal cone of nonnegative functions. It
is easy to verify that 
$$P^{\circ}=\{x\in P: \exists \epsilon>0 \hbox{ such that } x(t)\geq \epsilon, \hbox{ for all }t\in \R \}.$$

\begin{Def}(Guo and Lakshmikantham \cite{GuoL}) 
	Let $(X,\le)$ be an ordered Banach space 
	and let $E\subset X$. An operator $\Phi: E\times E \to X$ is called a \textbf{ mixed monotone operator } if $\Phi(x,y)$ is nondecreasing in $x$ and 
	nonincreasing 
	in $y$. 
	An element $\tilde{x}\in E$ is called a \textbf{fixed point } of $\Phi$ if $\Phi(\tilde{x},\tilde{x})=\tilde{x}.$
\end{Def}

We first establish the following abstract fixed point Lemma which will play an important role in sequel.

\begin{Lem}
	\label{puntofijolema16}
	Let $\Phi$ be an operator $\Phi:P^{\circ}\to P^{\circ}$. 
	Assume that 
	\begin{itemize}
		\item[$(I)$] there exist $u_0, v_0 \in P^{\circ}$ , $u_0 < v_0$  such that $u_0\le \Phi(u_0, v_0)$ and $v_0\geq \Phi(v_0,u_0)$;
		\item[$(II)$] $\Phi$ is a mixed monotone operator on $[u_0,v_0]\times [u_0, v_0]$;
		\item[$(III)$] there exists a function $\phi:\left[\frac{u_0^-}{ v_0^+} ,1\right)\to (0,+\infty)$ such that $\phi(\gamma)>\gamma$, for any $x, y \in [u_0,v_0]$ 
		
		$$\Phi(\gamma x, \gamma^{-1} y)\geq \phi(\gamma)\Phi(x, y), \,\,\, \hbox{ for all  } \gamma\in \left[{\frac{u_0^-}{ v_0^+}},1\right).$$
	\end{itemize}
	Then $\Phi$ has exactly one 
	fixed point $\tilde{x}$ in $[u_0,v_0]$.
\end{Lem}

\begin{proof}
	Let $u_n:=\Phi(u_{n-1}, v_{n-1})$ and $v_n:=\Phi(v_{n-1},u_{n-1})$, $n\in \N$. By   $(I)$ and the mixed monotonicity of $\Phi$, we have
	$u_0\le u_1= \Phi(u_0, v_0)\le \Phi(v_0, u_0) = v_1 \le v_0,$	
	and inductively  we obtain
	\begin{equation}
	\label{inequality1}
	u_0\le u_1\le \ldots \le u_n \le \ldots \le v_n \le \ldots \le v_1 \le v_0.
	\end{equation}
	Since $u_n$ is in the open set $P^{\circ}$,  
	there exists a constant 
	$\delta>0$ such that for any $\lambda\in (0,\delta)$, $u_n-\lambda v_n\in P^{\circ}$.
	
	Thus, the constant $\lambda_n:=\sup\{\lambda: u_n \geq \lambda v_n\}$ is well 
	defined for each $n\in \N_0$ and
	\begin{equation}
	\label{lambda_n1}
	u_n\geq \lambda_n v_n.
	\end{equation} 
	Moreover, from $u_{n+1}\geq u_n \geq \lambda_n v_n\ge \lambda_n v_{n+1}$,
	it is easy to verify that $\lambda_{n+1}\geq \lambda_n$. Thus, inductively we obtain
	$$\lambda_0\le \lambda_1 \le \cdots \le \lambda_n \le \cdots \le 1 .$$
	In addition, since $\frac{u_0^-}{v_0^+}\in \{\lambda: u_0 \geq \lambda v_0\}$, we have
	$ \frac{u_0^-}{v_0^+}\le \lambda_n\le 1$ for all $n\in \N_0$. 
	
	We claim that $\overline{\lambda}:= \lim_{n\to +\infty}\lambda_n = 1$. Otherwise, suppose that 
	$\overline{\lambda} < 1$ and consider the following two cases:
	\begin{itemize}
		\item[ Case $1$.] Suppose that there exists $\overline{n}$ such that $\lambda_{\overline{n}} = 
		\overline{\lambda}$. Then  $u_n \geq \overline{\lambda} v_n$ for all $n> \overline{n}$
		which, together with $(II), (III)$ and (\ref{inequality1}), yield
		$$u_{n+1}= \Phi(u_n, v_n)\geq \Phi(\overline{\lambda} v_n, \overline{\lambda}^{-1}u_n)\geq 
		\phi(\overline{\lambda}) \Phi(v_n, u_n)= \phi(\overline{\lambda}) v_{n+1}.$$
		Thus, $\lambda_{n+1}\geq \phi(\overline{\lambda}) > \overline{\lambda}$, which contradicts the fact that 
		$\lambda_{n+1}=\overline{\lambda}$.
		
		\noindent
		\item[ Case $2$.] Suppose that $\lambda_n < \overline{\lambda}$, for all $n$. In view of $(II), (III)$ and (\ref{inequality1}), 
		we have
		\begin{align*}		
		u_{n+1} &\geq \Phi\left(\lambda_n v_n, \lambda_n^{-1}u_n\right)=
		\Phi\left( \frac{\lambda_n}{\overline{\lambda}}\overline{\lambda} v_n, \frac{\overline{\lambda}}{\lambda_n}\overline{\lambda}^{-1}u_n\right)\\ 
		&\geq \phi\left(\frac{\lambda_n}{\overline{\lambda}}\right)\Phi\left(\overline{\lambda}v_n,\overline{\lambda}^{-1}u_n\right)>\frac{\lambda_n}{\overline{\lambda}}
		\phi\left(\overline{\lambda}\right) \Phi(v_n,u_n)\geq \frac{\lambda_n}{\overline{\lambda}}\phi\left(\overline{\lambda}\right)  v_{n+1}.
		\end{align*}

		Thus, $\lambda_{n+1}\geq \frac{\lambda_n}{\overline{\lambda}} \phi\left(\overline{\lambda}\right) $.
		Letting $n\to \infty$, we obtain $\overline{\lambda}\geq \phi\left(\overline{\lambda}\right) > \overline{\lambda}$, a contradiction. 
	\end{itemize}
	In view of 	(\ref{inequality1})-(\ref{lambda_n1}), it follows that 
	\begin{equation}
	\label{inequality2}
	0\le u_{n+k}-u_{n}\le v_n-u_n\le v_n-\lambda_n v_n = (1-\lambda_n)v_n\le (1-\lambda_n) v_0, \,\,\, \hbox{ for any } n,k\in\N.
	\end{equation}
	It is followed by the normality of $P$ and (\ref{inequality2}), 
	$$||u_{n+k}-u_n||\le N (1-\lambda_n) ||v_0||\to 0 \,\,\, \hbox{ as } n\to \infty.$$ 
	Thus, $\{u_n\}_{n\in \N}$ is a Cauchy sequence. This implies that 
	there exists $\overline{u}\in [u_0, v_0]$ such that $u_n\to \overline{u}$. Similarly,
	$$	0\le v_{n}-u_{n}\le v_n-\lambda_n v_n = (1-\lambda_n)v_n\le (1-\lambda_n) v_0$$
	and
	$$||v_n-u_n||\le N (1-\lambda_n) ||v_0||\to 0 \,\,\, \hbox{ as } n\to \infty,$$ 
	which implies that $v_n\to \overline{u}$.
	From the mixed monotonicity of $\Phi$ on $[u_0,v_0]$, we have
	$$u_{n+1}=\Phi(u_n, v_n) \le \Phi(\overline{u}, \overline{u})\le  \Phi(v_n, u_n) = v_{n+1.} $$
	We conclude that $\overline{u}=\Phi(\overline{u},\overline{u})$.
	
	Suppose now that $\overline{w}\in [u_0, v_0]$ is a fixed point of $\Phi$. Let
	$\alpha :=\sup\{\tilde{\alpha}\in(0,1):\tilde{\alpha} \overline{w} \le \overline{u} \le \tilde{\alpha}^{-1} \overline{w}\}$. Since 
	$\overline{u},\overline{w}$ have positive infimum, $\alpha$ is well defined.
	In addition, $\alpha \overline{w}\le \overline{u} \le \alpha^{-1} \overline{w}$
	and $\alpha \in \left[\frac{u_0^-}{v_0^+},1\right]$. Suppose that $\alpha \in \left[\frac{u_0^-}{v_0^+},1\right)$, then
	
	$$\overline{u}=\Phi(\overline{u},\overline{ u})\le \Phi(\frac{1}{\alpha} \overline{w},\alpha\overline{ w})\le \phi(\alpha)^{-1} \Phi(\overline{w},\overline{w})=\phi(\alpha)^{-1} \overline{w},$$
	and 
	$$\overline{u}=\Phi(\overline{u},\overline{ u})\geq \Phi(\alpha \overline{w},\frac{1}{\alpha}\overline{ w})\geq \phi(\alpha) \Phi(\overline{w},\overline{w})=\phi(\alpha) \overline{w}.$$
	Thus, by the definition of $\alpha$ we have $\phi(\alpha)\le \alpha$,  a contradiction. Therefore $\overline{w}=\overline{u}$, and the proof is complete.	
\end{proof}

\begin{Lem}
	\label{formulasol}
	Let	$x(t)$ be  a solution of 
	\begin{equation}
	\label{fun}
	x'(t) = f(t,x(t-\tau_1(t)),\cdots,x(t-\tau_M(t)))-b(t)x(t).
	\end{equation}
	If $x(t)$ is defined in the whole real line, then
	\begin{equation}
	\label{eq11}
	x(t) = \int_{-\infty}^{t}e^{-\int_{s}^t b(u) du}f(t,x(t-\tau_1(t)),\cdots,x(t-\tau_M(t)))ds , \,\,\, \hbox{ for all } t\in \R.
	\end{equation}
\end{Lem}

\begin{proof}
	Let $t_1,\in \R$. Integrating (\ref{fun}) from $t_1$ to $t$, we have
	$$
	x(t)  = x(t_1) e^{-\int_{t_1}^{t} b(u)du} + \int_{t_1}^t e^{-\int_s^t b(u) du}f(t,x(t-\tau_1(t)),\cdots,x(t-\tau_M(t))) ds.
	$$
	In addition, as $x(t)$ is defined on the whole real line, by taking the limit on the right-hand side of the equality we obtain (\ref{eq11}).
	The proof is  complete.
\end{proof}

\begin{Rem}
	Let $A, m, n$ be  positive constants, with $m>1$. Define the constant $B$ as
	\begin{equation}
	\label{B}
	B := A \left[\frac{n A^n}{(m-1)(1+A^n)}\right]^{\frac{1}{n-m+1}}.
	\end{equation}
	Thus,
	\begin{equation}
	\label{AB}
	B>A \hbox{ if and only if } A> \left(\frac{m-1}{n-m+1}\right)^{\frac{1}{n}}.
	\end{equation}
\end{Rem}

\section{Existence of almost periodic solutions.}
\label{main}

\begin{Th}
	\label{m>1}
	Let $n>m-1>0$ and $n\le m $. Let $A$ be a constant such that $A> \left(\frac{m-1}{n-m+1}\right)^{\frac{1}{n}}$ and $B$ defined in
	(\ref{B}). Furthermore, assume that
	\begin{equation}
	\label{cond}
	\frac{1+A^n}{A^{m-1}}\le \sum_{k=1}^M\frac{\lambda_k r_k^-}{b^+}\le F^s\sum_{k=1}^M \frac{\lambda_k r_k^+}{(b^*)^-}\le \frac{1+B^n}{B^{m-1}}.
	\end{equation}
	Then (\ref{MG}) has a unique almost periodic solution  $x(t)$ such that $A\le x(t)\le B$.
\end{Th}

\begin{proof} 
	Consider the positive constant functions $u_0 = A $
	and $v_0= B$.	Set the operator
	\begin{equation}
	\label{operator}
	\Phi(x)(t):= \int_{-\infty}^t e^{-\int_s^t b(u)du}\sum_{k=1}^M\lambda_kr_k(s) \frac{x^m(s-\tau_k(s))}{1+x^n(s-\tau_k(s))}ds\,\,\,\hbox{ for all } t \in \R.
	\end{equation}
	Clearly the function $f(u)=\frac{u^{m}}{1+u^{n}}$ is nondecreasing on  $[A,B]$. Thus, we deduce that $\Phi$ is a nondecreasing operator on $[u_0,v_0]$. It follows from  properties of almost periodic functions that $\Phi(P^\circ)\subset AP(\R)$.  Moreover,
	$\Phi$ satisfies $\Phi( P^{\circ}) \subset P^{\circ}$. Indeed,  
	let $x\in P^{\circ}$,  there exists $\epsilon >0$ such that $x(t) \geq \epsilon$
	for all $t\in \R$. Then,

	\begin{align*}
	\Phi(x)(t) & \geq \int_{-\infty}^t e^{-b^+(t-s)}\sum_{k=1}^M\lambda_k r_k^-
	\frac{\epsilon^{m}}{1+\epsilon^{n}}ds \geq \sum_{k=1}^M\frac{\lambda_k r_k^-}{b^+} \frac{\epsilon^{m}}{1+\epsilon^{n}} := \tilde{\epsilon}>0.
	\end{align*}
	
	Now, by virtue of (\ref{oscilatoria}) and (\ref{cond}) we find 
	
	\begin{align*}
	\Phi(u_0) &=\int_{-\infty}^t e^{-\int_s^t b(u)du}\sum_{k=1}^M\lambda_kr_k(s) \frac{u_0^m}{1+u_0^n}ds \geq \sum_{k=1}^M\frac{\lambda_k r_k^-}{b^+}\frac{u_0^{m}}{(1+u_0^{n})}\geq u_0,
	\end{align*}
	
	\noindent
	and	
	\begin{align*}
	\Phi(v_0)(t) 
	&\le F^s\sum_{k=1}^M\frac{\lambda_k r_k^+}{(b^*)^-}\frac{v_0^{m}}{1+v_0^{n}}
	\le v_0.
	\end{align*}

	Finally, it only remains to show that condition $(IV)$ of Lemma \ref{puntofijolema16} is satisfied.
	Moreover, for each $\gamma\in \left[\frac{A}{B},1\right)$, $x\in[u_0,v_0]$ and $t\in \R$, from the monotonicity of 
	$f(u)=\frac{1+u^n}{1+\gamma^n u^n}$, we have
	\begin{align*}
	\Phi(\gamma x)(t) & = \int_{-\infty}^t e^{-\int_s^t b(u)du}\sum_{k=1}^M\lambda_kr_k(s) \frac{x^m(s-\tau_k(s))}
	{1+ x^n(s-\tau_k(s))}\gamma^m\frac{1+ x^n(s-\tau_k(s))}{1+\gamma^n x^n(s-\tau_k(s))}ds\\
	& \geq \Phi(x)(t) \gamma^m\frac{1+ u_0^n}{1+\gamma^n u_0^n} := \Phi(x)(t) \phi(\gamma),
	\end{align*}

	\noindent
	where $\phi:\left[\frac{A}{B},1\right)\to (0,+\infty)$  is the mapping defined by
	\begin{equation}
	\label{functionphi}
	\phi(\gamma)= \gamma^m\frac{1+ A^n}{1+\gamma^n A^n}.
	\end{equation}
	Thus, $$\Phi(\gamma x)\geq \phi(\gamma)\Phi(x), \hbox{ for each } \gamma\in \left[\frac{A}{B},1\right) \hbox{ and } x\in [u_0,v_0].$$

	In order to prove that $\phi(\gamma)>\gamma$, for convenience, we define the function 
	$M(\gamma):=\gamma^{m-1}(1+ A^n)-(1+\gamma^n A^n).$
	It is easy to verify  that $M(\gamma)$ achieves the maximum in $\gamma_{\max}=\left[\frac{(m-1)(1+A^n)}{n A^n}\right]^{\frac{1}{n-m+1}}=\frac{A}{B}$,  $M(1)=0$ and 
	$M(\gamma)$ is strictly decreasing in 
	$\left(\frac{A}{B},1\right)$, which implies that $M(\gamma)>0 $ for all $\gamma \in \left[\frac{A}{B},1\right)$. Thus,  $\phi(\gamma)>\gamma$ for $\gamma \in \left[\frac{A}{B},1\right)$.
	Thus, $\Phi$ satisfies all the assumptions of Lemma \ref{puntofijolema16} and $\Phi$ has a unique fixed point $x\in [u_0,v_0]\subset P^{\circ}$. By Lemma \ref{formulasol}, $x(t)$ is the unique almost periodic solution of (\ref{MG}) which satisfies $A\le x(t)\le B.$
\end{proof}

\begin{Th}
	\label{m>12}
	Let $n > m> 1 $. Let $A$ be a constant such that $ \left(\frac{m-1}{n-m+1}\right)^{\frac{1}{n}}< A \le \left(\frac{m}{n-m}\right)^{\frac{1}{n}}\le B,$ with $B$ defined in
	(\ref{B}).  Assume that
	\begin{equation}
	\label{cond1}
	\frac{1+A^n}{A^{m-1}}\le \sum_{k=1}^M \frac{\lambda_k r_k^-}{b^+}\le F^s \sum_{k=1}^M \frac{\lambda_k r_k^+}{(b^*)^-}\le \left(\frac{m}{n-m}\right)^{\frac{1}{n}}.
	\end{equation}
	Then (\ref{MG}) has a unique almost periodic solution  $x(t)$ such that $x(t)\geq A$. 
\end{Th}


\begin{proof}
	The proof is divided into 2 steps.
	
	\noindent
	\textit{Step 1.} 	
	Let $x(t)$ an almost periodic solution of (\ref{MG}). According to Lemma \ref{formulasol}, (\ref{oscilatoria}), (\ref{cond1}) and the fact that
	\begin{equation}
	\label{r1}
	\sup_{u\geq 0}\left\{\frac{u^{m}}{1+u^{n}}\right\}\le 1, \hbox{ for } n\geq m\geq 0,
	\end{equation}	
	we have
	
	\begin{align*}
	x(t) & \le \int_{-\infty}^t e^{-\int_s^t b(u) du}\sum_{k=1}^M \lambda_k r_k(s) ds \le F^s\sum_{k=1}^M \frac{\lambda_k r^+_k}{(b^*)^-} \le \left(\frac{m}{n-m}\right)^{\frac{1}{n}}:=V, \qquad \text{ for all } t\in \R.
	\end{align*}

	\noindent
	\textit{Step 2.} 
	Let us define the following truncated function $h$ for $x>0$, namely 
	
	\begin{equation}
	\label{func}
	h(x):= \left\{ \begin{array}{lcc}
	\frac{x^{m}}{1+x^{n}} &   \hbox{ if } & x \leq V \\
	\\ \frac{V^{m}}{1+V^{n}} &  \hbox{ if } &   x > V 
	\end{array}
	\right.,
	\end{equation}
	with $V$ defined in \textit{Step 1.}
	Let us consider the following associated equation:
	\begin{equation}
	\label{integral2}
	x'(t) = \sum_{k= 1}^M \lambda_k r_k(t)h(x(t-\tau_k(t)))-b(t) x(t).
	\end{equation}
	We define the nonlinear operator $\Theta$ on $P^{\circ}$ by,
	$$\Theta(x)(t):=\int_{-\infty}^t e^{-\int_{s}^t b(u) du}\sum_{k=1}^M \lambda_k r_k(s)h(x(s-\tau_k(s)))ds, \,\,\, t\in \R.$$
	
	Let $u_0:= A$ and $v_0:= B$. It is not difficult to prove that $\Theta $ is a non-decreasing operator and $\Theta(P^{\circ}\times P^{\circ})\subset P^{\circ}.$
	
	For each $x\in[A,B]$ and $\gamma\in \left[\frac{A}{B},1\right)$, we have

	\begin{equation}
	\label{f}
	\frac{h(\gamma x)}{h(x)}\geq\left\{\begin{array}{lcc}
	\phi(\gamma) &   \hbox{ if } & x \leq\frac{1}{\gamma} V \\
	\\1 &  \hbox{ if } &   x >\frac{1}{\gamma} V,
	\end{array}
	\right.
	\end{equation}
	where $\phi$ is the same as in (\ref{functionphi}). Letting $\theta(\gamma) =\min\{\phi(\gamma),1 \}$,  it is readily verified that for each $\gamma \in \left[\frac{A}{B},1\right)$ and $x\in [u_0,v_0]$,
	$$\Theta(\gamma x)\geq \theta(\gamma)\Theta(x),$$
	where $\theta(\gamma):=\min\{\phi(\gamma),1 \}$. Analogously to the proof in Theorem \ref{m>1}, we can show that $\theta(\gamma)>\gamma$ for all $\gamma \in \left[\frac{A}{B},1\right)$, and that the remaining assumptions of Lemma \ref{puntofijolema16} are satisfied. Thus, $\Theta$ has a unique fixed point $\tilde{x}\in [u_0,v_0]$. In addition, from Lemma \ref{formulasol} and (\ref{cond1}), we get
	\begin{align*}
	\tilde{x}(t) & = \int_{-\infty}^t e^{-\int_s^t b(u)du}\sum_{k=1}^M \lambda_kr_k(s)h(\tilde{x}(s-\tau_k(s)))ds \le F^s \sum_{k=1}^M \frac{\lambda_k r_k^+}{(b^*)^-}\le V 
	\end{align*}

	\noindent
	which yields that $h(\tilde{x}(s-\tau_k(s)))= \frac{\tilde{x}^m(s-\tau_k(s))}{1+\tilde{x}^n(s-\tau_k(s))}$. Thus, again by Lemma \ref{formulasol}, $\tilde{x}$ is a solution for $(\ref{MG})$, with $\tilde x(t) \geq A$ for all $t\in \R$.
	
	Let $\tilde z$ be an almost periodic solution for (\ref{MG}) such that 
	$\tilde z(t)\geq A$ for all $t\in \R$. Thus, by \textit{Step 1} we conclude that $\tilde z$ is an almost periodic solution for  (\ref{integral2}) such that $A\le \tilde z \le V\le B$, which means that $\tilde x = \tilde z$. The proof is complete.
\end{proof}

\medskip
\noindent

\section{Examples}
Consider the following model of hematopoiesis with multiple time-varying delays:

\begin{equation}
\label{ejemplo}
x'(t) = \frac{1}{2}\left(5+|\cos(\sqrt{2}t)|\right)\frac{x^{m}\left(t-2e^{\cos t}\right)}{1+x^{n}\left(t-2e^{\cos t}\right)} +\frac{1}{4}\left(13+\frac{3}{5}|\sin(\sqrt{3}t)|\right)\frac{x^{m}\left(t-2e^{\sin t}\right)}{1+x^{n}\left(t-2e^{\sin t}\right)}-(1+1.2\cos(400t)) x(t).
\end{equation}

It is seen that, 
$ b(t) = 1+1.2\cos(400t), \,\,\,M[b] = 1,\,\,\, b^+ = 2.2, \,\,\, b^*(t) = 1,\,\,\, (b^*)^- = 1 \hbox{ and } F^s = e^{\frac{1.2}{200}}.$ 
\begin{Ex}	Consider $ m = \frac{11}{10}$ and $n= \frac{1}{2}$ in (\ref{ejemplo}). Let $A := 4 > \left(\frac{m-1}{n-m+1}\right)^{\frac{1}{n}} = 0.0625,$ so $B = 4\left(\frac{10}{3}\right)^{\frac{5}{2}} \approx 81.14408$ and
	$$ 2.6116 \approx\frac{1+A^n}{A^m-1} <  \frac{5.75}{2.2} = \frac{\sum_{k=1}^M \lambda_k r^-_k}{b^+}\le F^s\frac{\sum_{k=1}^M \lambda_k r^+_k}{(b^*)^-}= e^{\frac{1.2}{200}}6.4  <  \frac{1+B^n}{B^m-1}\approx 6.4479 .$$
	Thus, all assumptions of Theorem \ref{m>1} are satisfied. Therefore, equation (\ref{ejemplo}) has a unique almost periodic solution $x(t)$, which satisfies $4\le x(t)\le 4\left(\frac{10}{3}\right)^{\frac{5}{2}}$. 
	
\end{Ex}

\begin{Ex}Consider $ m = \frac{11}{10}$ and $n= \frac{12}{10}$ in (\ref{ejemplo}). Let $A := 1.3>\left(\frac{m-1}{n-m+1}\right)^{\frac{1}{n}} \approx 0.13557,$ so $B \approx 7.56193$ and
	$$\frac{1+A^n}{A^m-1}\approx 2.30866\le 5.75 = \frac{\sum_{k=1}^M \lambda_k r^-_k}{b^+}\le F^s\frac{\sum_{k=1}^M \lambda_k r^+_k}{(b^*)^-}= F^s6.4 \approx 6.4385\le  \left(\frac{m}{n-m}\right)^{\frac{1}{n}}\le B.$$
	Thus, all assumptions of Theorem \ref{m>12} are  satisfied. Therefore, equation (\ref{ejemplo}) has a unique almost periodic solution $x(t)$ such that $ x(t)\geq 1.3$. 
\end{Ex}

\section{Concluding remarks and open problem}

The fixed point theorems,  employed in previous works to establish existence results for the hematopoiesis model, involve functions such as $\phi:(0,1)\to (0,+\infty)$ or 
$\phi:(0,1)\times P^{\circ}\times P^{\circ}\to (0,+\infty)$ 
satisfying $\phi(\gamma)>\gamma$ or $\phi(\gamma,x,y)>0$ for all $x,y \in P^{\circ} $ and $\gamma\in(0,1)$ (see for example \cite{ANS, DLN, DZ, WZ} ). Unfortunately, these theorems fail when $m>1$  since the aforemention assumptions about $\phi$ are not fulfilled when $\gamma \approx 0$. Moreover, to the best of our knowledge, the authors only study the model of hematopoiesis with $m\le 1$. This implies that our results are new and they complement previously known results.

It is worth to notice that the approach used in this paper cannot be applied to equation (\ref{MG}) with $0< n\le m-1$. In fact, for $u_0,v_0\in P^{\circ}$, such that $u_0 < v_0$  conditions
$\Phi(u_0)\geq u_0, \,\,\,\Phi(v_0)\le v_0$ cannot be satisfied simultenously, as it is required in Lemma \ref{puntofijolema16}. It is an open and interesting problem.






\end{document}